\documentclass[11pt,a4paper]{article}

\usepackage[a4paper,margin=30mm]{geometry}
\usepackage[T1]{fontenc}
\usepackage{lmodern}
\usepackage{microtype}

\usepackage[reqno]{amsmath}
\usepackage{amssymb,amsthm,mathtools,comment}
\usepackage{aliascnt}

\usepackage{enumitem}
\usepackage[numbers,sort&compress]{natbib}
\usepackage{tikz}
\usetikzlibrary{arrows.meta}
\usepackage[colorlinks=true,linkcolor=blue,citecolor=blue,urlcolor=blue]{hyperref}
\usepackage[noabbrev]{cleveref}
\makeatletter
\renewcommand{\section}{\@startsection{section}{1}\z@
  {2.5ex\@plus 1ex\@minus .2ex}{1.2ex\@plus .2ex}
  {\normalfont\large\bfseries\centering}}
\renewcommand{\subsection}{\@startsection{subsection}{2}\z@
  {2ex\@plus .8ex\@minus .2ex}{.8ex\@plus .2ex}
  {\normalfont\normalsize\bfseries\raggedright}}
\makeatother

\renewenvironment{abstract}{%
  \par\medskip
  {\centering\normalfont\normalsize\bfseries Abstract\par}%
  \vspace{0.6em}%
  \begingroup
  \normalfont\normalsize
  \leftskip=3pc\relax
  \rightskip=3pc\relax
  \parindent=1.5em\relax
  \indent\ignorespaces
}{%
  \par\endgroup\medskip
}

\newtheorem{theorem}{Theorem}[section]
\newaliascnt{lemma}{theorem}
\newtheorem{lemma}[lemma]{Lemma}
\aliascntresetthe{lemma}
\newaliascnt{proposition}{theorem}
\newtheorem{proposition}[proposition]{Proposition}
\aliascntresetthe{proposition}
\newaliascnt{corollary}{theorem}

\aliascntresetthe{corollary}
\newaliascnt{claim}{theorem}
\newtheorem{claim}[claim]{Claim}
\aliascntresetthe{claim}
\newaliascnt{conjecture}{theorem}
\newtheorem{conjecture}[conjecture]{Conjecture}
\aliascntresetthe{conjecture}
\newaliascnt{definition}{theorem}
\newtheorem{definition}[definition]{Definition}
\aliascntresetthe{definition}
\newaliascnt{fact}{theorem}
\newtheorem{fact}[fact]{Fact}
\aliascntresetthe{fact}
\theoremstyle{remark}
\newaliascnt{remark}{theorem}

\aliascntresetthe{remark}

\crefname{theorem}{Theorem}{Theorems}
\crefname{lemma}{Lemma}{Lemmas}
\crefname{claim}{Claim}{Claims}
\crefname{section}{Section}{Sections}
\crefname{corollary}{Corollary}{Corollaries}
\crefname{proposition}{Proposition}{Propositions}
\crefname{conjecture}{Conjecture}{Conjectures}
\crefname{fact}{Fact}{Facts}
\crefname{equation}{Equation}{Equations}
\crefname{definition}{Definition}{Definitions}
\crefname{remark}{Remark}{Remarks}

\newcommand{\eps}{\varepsilon}
\newcommand{\dplus}{d^{+}}
\newcommand{\dminus}{d^{-}}
\newcommand{\dsemi}{\delta^{0}}
\newcommand{\RNp}{RN^{+}}
\newcommand{\dist}{\operatorname{dist}}

\newcommand{\cG}{\mathcal G}

\newcommand{\dmax}{\Delta^{0}}

\allowdisplaybreaks

\title{Cycles in Regular Multipartite Tournaments}
\author{%
  Yufei Chang$^{\dagger}$\quad
  Zhilan Wang$^{\dagger}$\quad
  Shuo Wei$^{\dagger}$\quad
  Jin Yan$^{\dagger}$
}
\date{}

\begin{document}
\maketitle

\begingroup
\renewcommand{\thefootnote}{\fnsymbol{footnote}}
\footnotetext[2]{School of Mathematics, Shandong University,
Jinan 250100, China. This work was supported by the National Natural Science
Foundation of China (No.~12571373).}
\footnotetext[1]{\fontsize{8.5pt}{10.5pt}\selectfont
\raggedright
 Email addresses:
\texttt{202420326@mail.sdu.edu.cn} (Yufei Chang);
\texttt{202690000074@sdu.edu.cn} (Zhilan Wang);
\texttt{shuowei@mail.sdu.edu.cn} (Shuo Wei);
\texttt{yanj@sdu.edu.cn} (Jin Yan).}
\endgroup

\begin{abstract}
We study balanced cycles in regular multipartite tournaments. A cycle in a $c$-partite tournament is $k$-balanced if it contains exactly $k$ vertices from each partite class. We prove that there exists an absolute integer $r_0$ such that every regular $c$-partite tournament with each partite class size of $r\ge r_0$ contains a $k$-balanced cycle for every $1\le k\le r$ when $c=3$ or $c\ge5$, and for every $2\le k\le r$ when $c=4$. In particular, our result confirms three conjectures of Volkmann and Winzen for large regular multipartite tournaments. We also confirm Yeo's conjecture on the existence of $(p,q)$-cycle factors for large regular $c$-partite tournaments $G$ with $c\ge4$, where a $(p,q)$-cycle factor is a pair of vertex-disjoint cycles of lengths $p$ and $q$ spanning $G$.
\end{abstract}


\section{Introduction}\label{sec:introduction}

Cycles in tournaments have been widely studied. Camion~\cite{Camion} proved that every strongly connected tournament has a Hamilton cycle. Moon~\cite{Moon} established a more general result, showing that every strongly connected tournament on $n$ vertices is vertex-pancyclic, meaning that every vertex lies on a cycle of each length from $3$ to $n$.

A \emph{$c$-partite tournament} is an orientation of a complete $c$-partite graph. For a vertex $v$, let $d^+(v)$ and $d^-(v)$ denote its \emph{outdegree} and \emph{indegree}, respectively. A $c$-partite tournament is $d$-\emph{regular} if $d^+(v)=d^-(v)=d$ for every vertex $v$. In this case, all partite classes have the same size $r$, and $d=(c-1)r/2$. A natural question is whether similar pancyclicity results hold for regular $c$-partite tournaments. Yeo~\cite{YeoVertex} proved that every regular $c$-partite tournament is vertex-pancyclic when $c\ge5$. However, this fails for $c=3$. Indeed, consider the regular $3$-partite tournament with all edges directed from $V_i$ to $V_{i+1}$, where indices are taken modulo $3$. Every cycle in this regular $3$-partite tournament has length divisible by three. This motivates studying the existence of cycles of all lengths divisible by three. Bondy~\cite{Bondy}, Yeo~\cite{YeoHamilton} and Volkmann~\cite{VolkmannComplementary} established the existence of cycles of lengths $3$, $3r$ and $3r-3$, respectively, in regular tripartite tournaments. Volkmann~\cite{Volkmann} further proved the existence of a  $6$-cycle and proposed the following conjecture.

\begin{conjecture}[Volkmann~\cite{Volkmann}]\label{vo}
Every regular tripartite tournament with $r$ vertices in each class contains a cycle of length $3k$ for every $1\le k\le r$.
\end{conjecture}

Stella and Volkmann~\cite{StellaVolkmannNine} proved the case $k=3$ later. 
A path or cycle is \emph{$k$-balanced} if it contains exactly $k$ vertices from each partite class, and simply \emph{balanced} if $k=1$. 
Volkmann and Winzen~\cite{VolkmannWinzenCycles}
proved that every regular $c$-partite tournament
with $c\ge5$ contains an $(r-1)$-balanced cycle.
More generally, they showed that, for every fixed
positive integer $s$, such a tournament contains
an $(r-s)$-balanced cycle whenever $r$ is
sufficiently large. They further conjectured
the following.

\begin{conjecture}[Volkmann and Winzen \cite{VolkmannWinzenCycles}]\label{vwc}
Every regular $c$-partite tournament with $c\ge5$ and common partite size $r$ contains a $k$-balanced cycle for every $1\le k\le r$.
\end{conjecture}

 The case $k=1$ was later
settled in~\cite{VolkmannWinzenStrong}.
The condition $c\ge5$ is necessary, as Volkmann
constructed a family of regular four-partite
tournaments containing no balanced
cycle~\cite{VolkmannWinzenStrong}.
This motivated Volkmann and Winzen~\cite{vp}
to propose the following weaker conjecture
for balanced paths.

\begin{conjecture}[Volkmann and Winzen~\cite{vp}]
\label{vw}
Every regular $c$-partite tournament with
$c\ge2$ and common partite size $r$ contains
a $k$-balanced path for every
$1\le k\le r$.
\end{conjecture}

The conjecture is known to hold for $c=2$,
and the cases $k\in\{1,2,r-1,r\}$ have been
established for all $c\ge2$ in \cite{vp}.
Our first result confirms \Cref{vo,vwc,vw} for all $r\ge r_0$, where $r_0$ is an absolute constant.

\begin{theorem}\label{thm:unified}
There exists a positive integer $r_0$ such
that for every $c\ge3$ and $r\ge r_0$,
every regular $c$-partite tournament with
common partite size $r$ contains a $k$-balanced cycle
for every $1\le k\le r$ if $c=3$ or $c\ge5$,
and for every $2\le k\le r$ if $c=4$.
\end{theorem}

We also study \emph{$(p,q)$-cycle-factors} in regular $c$-partite tournaments. A $(p,q)$-cycle-factor consists of two vertex-disjoint cycles of lengths $p$ and $q$ that together cover all vertices. Yeo~\cite{YeoVertex} proposed the following conjecture.

\begin{conjecture}[Yeo~\cite{YeoVertex}]
\label{conj:yeo-factor}
Every regular $c$-partite tournament $G$
with $c\ge5$ contains a
$(p,|V(G)|-p)$-cycle-factor for every
$3\le p\le |V(G)|-3$.
\end{conjecture}

Bessy and Thiebaut~\cite{BessyThiebaut2023}
proposed a weaker conjecture under a stronger
regularity condition.

\begin{conjecture}[Bessy and Thiebaut~
\cite{BessyThiebaut2023}]
\label{conj:bessy-factor}
Every $c$-partite tournament $G$ with $c\ge5$,
in which the bipartite tournament between
any two partite classes is $t$-regular,
contains a $(p,|V(G)|-p)$-cycle-factor
for every even integer
$4\le p\le |V(G)|-4$.
\end{conjecture}

Clearly, \Cref{conj:yeo-factor} implies
\Cref{conj:bessy-factor}.
Volkmann~\cite{VolkmannFour2004} also
investigated the four-partite case,
particularly cycle-factors containing a $4$-cycle.
Our second main result confirms Yeo's
conjecture for sufficiently large tournaments
and extends it to $c=4$.

\begin{theorem}\label{thm:complementary}
There exists a positive integer $r_0$ such
that for every $c\ge4$ and $r\ge r_0$,
every regular $c$-partite tournament $G$ with
common partite size $r$ contains a $(p,cr-p)$-cycle-factor for every $3\le p\le cr-3$.
\end{theorem}

\medskip
\noindent\textbf{Organization.}
In \Cref{sec:prep}, we collect the preliminary results. In \Cref{sec:multipartite}, we prove \Cref{thm:unified} for $c\ge4$.
In \Cref{sec:tripartite}, we establish the tripartite case by applying the stability theorem \Cref{thm:stability},
with the robust outexpander and extremal cases handled in \Cref{lem:robust-case,lem:extremal-transfer}, respectively.
Finally, in \Cref{sec:complementary}, we prove
\Cref{thm:complementary}.

\section{Preliminaries}\label{sec:prep}

\subsection{Notations}\label{subsec:notation}
 Let $G$ be a digraph.
For $v\in V(G)$, let
$N_G^+(v)=\{u:vu\in E(G)\}$ and
$N_G^-(v)=\{u:uv\in E(G)\}$
be its out- and in-neighborhoods, respectively, and
$d_G^\pm(v)=|N_G^\pm(v)|$.
For $B\subseteq V(G)$, set
$d_G^\pm(v,B)=|N_G^\pm(v)\cap B|$ and
$d_G(v,B)=d_G^+(v,B)+d_G^-(v,B)$.
The total degree of $v$ is
$d_G(v)=d_G^+(v)+d_G^-(v)$, and the minimum
total degree of $G$ is $\delta(G)$.
The minimum and maximum semidegrees are $\dsemi(G)=\min_{v\in V(G)}
\min\{d_G^+(v),d_G^-(v)\}$ and $\dmax(G)=\max_{v\in V(G)}
\max\{d_G^+(v),d_G^-(v)\}.$
A vertex $v$ of a tournament $T$ is a
\emph{source} if $d_T^-(v)=0$ and a
\emph{sink} if $d_T^+(v)=0$.

For $x,y\in V(G)$, write $x\to y$ (or $xy$)
if $xy\in E(G)$, and $x\to y\to z$
if $xy,yz\in E(G)$.
For disjoint sets $A,B\subseteq V(G)$,
write $A\to B$ if $a\to b$ for all
$a\in A$ and $b\in B$.
For any $A,B\subseteq V(G)$, let $E_G(A,B)=\{xy\in E(G):x\in A,\ y\in B\},$
and
$e_G(A,B)=|E_G(A,B)|
=\sum_{v\in A}d_G^+(v,B)$.
We omit the subscript $G$ when clear.
An $(A,B)$-path (resp.\ $(A,B)$-edge)
is a path (resp.\ edge) from $A$ to $B$.
In particular, an $(a,b)$-path is a path from
$a$ to $b$.

A \emph{linear forest}
$\mathcal P=P_1\cup\cdots\cup P_t$
consists of nontrivial vertex-disjoint paths,
where each $P_i$ is an $(a_i,b_i)$-path.
Let $\mathring V(\mathcal P)$ be the set
of its internal vertices. Then $|V(\mathcal P)|\le 2e(\mathcal P)$.
We write $G\circ\mathcal P$ for the digraph
obtained by deleting $\mathring V(\mathcal P)$,
adding the edges $a_ib_i$ for all $i\in[t]$,
and retaining all other edges between
the remaining vertices.
These added edges form the \emph{auxiliary matching}
$M_{\mathcal P}:=\{a_ib_i:i\in[t]\}$.

For disjoint sets $A,B$, let $M[A,B]$
denote a matching directed from $A$ to $B$
with $|M[A,B]|=|A|=|B|$.

A hypergraph $\mathcal H$ consists of a vertex
set $V(\mathcal H)$ and an edge set
$E(\mathcal H)\subseteq 2^{V(\mathcal H)}$.
It is \emph{$k$-uniform} if every edge has size
$k$, and \emph{$c$-partite} if its vertices
can be partitioned into $c$ classes,
each meeting every edge in at most one vertex.
The complete $c$-partite $k$-uniform
hypergraph $K^{(k)}_{\ell_1,\ldots,\ell_c}$
has vertex classes of sizes
$\ell_1,\ldots,\ell_c$ and contains all
$k$-sets meeting each class in at most
one vertex.

Throughout the paper, we omit all floor and ceiling functions whenever these are not crucial. The constants in the hierarchies used to state our results are chosen from right to left. For example, if we claim that a result holds whenever $ 0 < a \ll b \ll c \leq 1 $, then there are non-decreasing functions $f,g:(0,1]\to(0,1]$ such that the result holds for all $ 0 < a, b, c \leq 1 $ with $ b \leq f(c) $ and $ a \leq g(b) $.

\subsection{Propositions in multipartite tournaments}
We collect some basic properties and results on multipartite tournaments.
\begin{proposition}[R\'edei \cite{Redei}]\label{thm:redei-tool}
Every tournament has a Hamilton path.
\end{proposition}

\begin{proposition}[Camion \cite{Camion}]\label{thm:camion-tool}
Every strongly connected tournament has a Hamilton cycle.
\end{proposition}

\begin{proposition}[Yeo \cite{YeoVertex}]\label{thm:yeo-vertex-tool}
Every regular $c$-partite tournament with $c\ge5$ is
vertex-pancyclic.
\end{proposition}

\begin{proposition}[Yeo \cite{Yeo2007}]
\label{thm:yeo-four-tool}
Every sufficiently large regular four-partite tournament is
vertex-pancyclic.
\end{proposition}

\begin{proposition}[Yeo \cite{YeoNear}]\label{thm:yeo-hamilton-tool}
Let $G$ be a $c$-partite tournament on $n$ vertices whose partite sizes satisfy $n_1\le\cdots\le n_c$.  If $\dmax(G)-\dsemi(G)\le(n-n_{c-1}-2n_c+2)/2,$
then $G$ contains a Hamilton cycle.
\end{proposition}

\begin{proposition}\label{reg}
Let $G$ be a regular $c$-partite tournament
with partite classes $V_1,\ldots,V_c$,
where $c\ge 2$. Then all partite classes
have the same size $r$.
Moreover, every vertex $v\in V(G)$ satisfies $d^+(v)=d^-(v)=(c-1)r/2.$
\end{proposition}

\subsection{Regularity and embedding lemmas}

For disjoint vertex sets $X,Y$ in a graph $G$, write
$d_G(X,Y)=e_G(X,Y)/(|X||Y|)$, where $e_G(X,Y)$ denotes the number of edges between $X$ and $Y$. The pair $(X,Y)$ is $\eps$-regular if
$|d_G(X',Y')-d_G(X,Y)|\le\eps$ whenever
$X'\subseteq X$, $Y'\subseteq Y$, $|X'|\ge\eps|X|$ and $|Y'|\ge\eps|Y|$. It is $(\eps,d)$-super-regular if it is $\eps$-regular, has density at least $d$, 
and furthermore, $d_G(a,Y)>d|Y|$ for all $a\in X$
and $d_G(b,X)>d|X|$ for all $b\in Y$.

\begin{lemma}[Degree form of the Regularity Lemma]
\label{lem:regularity}
For every integer $c\ge2$, $\varepsilon\in(0,1)$ and
$M'\in\mathbb{N}$, there exist positive integers
$M=M(\varepsilon,M',c)$ and $n_0=n_0(\varepsilon,M',c)$
such that the following holds.
Let $G$ be a $c$-partite graph with partite classes
$V_1,\ldots,V_c$ of common size $r$, where $n=cr\ge n_0$,
and let $d\in[0,1]$.
Then there exist an integer $t$ with $M'\le t\le M$,
a spanning subgraph $G'\subseteq G$, and a partition $V(G)=V_0\cup
 \bigcup_{i=1}^{c}\bigcup_{j=1}^{t}V_{ij},$ where $V_{ij}\subseteq V_i,$
such that the following hold.

\medskip

\textnormal{($i$)}
$|V_0\cap V_i|\le\varepsilon r$ for every $1\le i\le c$.

\textnormal{($ii$)}
All clusters $V_{ij}$ have a common size $m$.

\textnormal{($iii$)}
For every $x\in V(G)$, $\deg_{G'}(x)>\deg_G(x)-(d+\varepsilon)n.$

\textnormal{($iv$)}
$G'[V_{ij}]$ is empty for every $1\le i\le c$ and
$1\le j\le t$.

\textnormal{($v$)}
For any two distinct clusters $V_{ij}$ and $V_{i'j'}$,
the bipartite graph $G'[V_{ij},V_{i'j'}]$ is
$\varepsilon$-regular and has density either $0$ or at least $d$.
\end{lemma}

The reduced graph $R$ has the clusters $V_{ij}$ as its vertices, with two distinct clusters adjacent whenever $G'$ contains an edge between them.

The following Blow-up Lemma was proved by Koml\'{o}s, S\'{a}rk\"{o}zy and Szemer\'{e}di \cite{KSS}.

\begin{lemma}[Blow-up Lemma~\cite{KSS}]\label{blowup}
For every $d>0$ and $\Delta\in\mathbb N$, there exists
$\tau_0=\tau_0(d,\Delta,F)>0$ such that the following holds.
Let $F$ be any graph and let $\ell\in\mathbb N$. Let $F(\ell)$
be the $\ell$-blow-up of $F$.
Suppose that $G_1$ is a spanning subgraph of $F(\ell)$ such that,
for every $xy\in E(F)$, the pair
$G_1[V_i,V_j]$ is $(\tau,d)$-super-regular, where
$0<\tau\le\tau_0$.
Then $G_1$ contains a copy of every spanning subgraph
$H\subseteq F(\ell)$ with $\Delta(H)\le\Delta$,
respecting the vertex classes.
\end{lemma}

Since our definition of super-regularity differs slightly from that in~\cite{KOT2005}, we make a minor adjustment to the parameters in the following lemma.

\begin{lemma}[K\"uhn, Osthus and Taraz \cite{KOT2005}]
\label{prop:super-regularization}
Let $G$ be a graph, and let $G'$ and $R$ be
the pure graph and the reduced graph,
respectively, obtained from the degree form
of the Regularity Lemma with parameters
$\varepsilon$ and $d$.
Suppose the nonexceptional clusters
$V_1,\ldots,V_t$ have common size $m$.
Let $H$ be a subgraph of $R$ satisfying
$\Delta(H)\le \Delta$, where
$\Delta\varepsilon<1$.
Then there exist subsets
$V_i'\subseteq V_i$ for every
$V_i\in V(H)$, each of size
$(1-\Delta\varepsilon)m$, such that
for every $V_iV_j\in E(H)$,
the pair $(V_i',V_j')$ in $G'$ is $(2\varepsilon/(1-\Delta\varepsilon),d-(\Delta+1)\varepsilon)$-super-regular.
\end{lemma}

\begin{lemma}[Erd\H{o}s~\cite{E}]
\label{thm:erdos-stone-multipartite}
For every $c,k\in\mathbb{N}$ with $c\ge2$
and every $\alpha>0$, there exists $r_0$
such that the following holds.
Let $H$ be a $c$-partite $c$-uniform
hypergraph with vertex classes
$V_1,\ldots,V_c$, each of size $r\ge r_0$.
If $e(H)\ge\alpha r^c$, then $H$ contains
a copy of $K^{(c)}_{k,\ldots,k}$.
\end{lemma}

\subsection{Robust outexpanders}\label{subsec:tripartite-stability}

For an $n$-vertex digraph $D$, a set $S\subseteq V(D)$ and $\nu>0$, define $\RNp_{\nu,D}(S):=\{v\in V(D):|N_D^-(v)\cap S|\ge \nu n\}.$
The digraph is a \emph{robust $(\nu,\tau)$-outexpander} if
$|\RNp_{\nu,D}(S)|\ge |S|+\nu n$ whenever
$\tau n\le |S|\le(1-\tau)n$. We use the following results on outexpanders to prove the non-extremal case.

\begin{lemma}[Gir\~ao, Granet, K\"uhn, Lo, and  Osthus \cite{GGKLO}]\label{lem:random-robust-expansion}
Let $0<1/n\ll\xi\ll\nu'
\ll\alpha,\nu,\tau\ll1$.
Suppose that $G$ is a robust
$(\nu,\tau)$-outexpander on $n$ vertices
and satisfies $\delta^0(G)\ge\alpha n$.
Then, for every integer $\xi n\le m\le n$,
a subset $S$ chosen uniformly at random
from all $m$-subsets of $V(G)$ satisfies $\mathbb P\bigl(
G[S]\text{ is a robust }
(\nu',4\tau)\text{-outexpander}
\bigr)\ge 1-1/n^2.$
\end{lemma}

\begin{lemma}[Lo and Patel \cite{LoPatel}]\label{thm:lo-patel}
Let $0<\nu\le\tau\le\gamma/16<1/16$. There exists a positive integer $n_0$ such that for all $n\ge n_0$,
every $n$-vertex digraph $D$ with $\dsemi(D)\ge\gamma n$ which is a robust
$(\nu,\tau)$-outexpander contains a cycle of every length $\nu n/2\le\ell\le n$.
\end{lemma}

\subsection{Other useful tools}

\begin{lemma}[Hoeffding inequality \cite{Hoe}]\label{cher}
Let $X=X_1+\cdots+X_q$, where the $X_i$
are independent random variables, each having
a binomial or hypergeometric distribution
with sample size $m_i$.
Set $M=\sum_{i=1}^q m_i$.
Then, for every $t\ge 0$, $\mathbb P\bigl(|X-\mathbb EX|\ge t\bigr)
\le 2\exp\left(-\frac{2t^2}{M}\right).$
\end{lemma}

\begin{lemma}[Markov's inequality]
Let $X$ be a nonnegative random variable. Then, for every $a>0$, $\mathbb{P}[X\ge a]
\le \mathbb{E}[X]/a.$
\end{lemma}

\begin{proposition}\label{log}
Let $r,t$ be positive integers and let $0\le k\le r$.
Then $\binom{tr}{tk}\le (r+1)^{t-1}\binom{r}{k}^{t}.$
\end{proposition}

\begin{proof}
Let $a_1,\ldots,a_t$ be integers such that $0\le a_i\le r$ for every $i\in[t]$ and $\sum_{i=1}^t a_i=tk$. Then $\prod_{i=1}^{t}\binom{r}{a_i}
\le \binom{r}{k}^{t}.$
Since there are at most $(r+1)^{t-1}$ choices for
$(a_1,\ldots,a_t)$, Vandermonde's identity gives
\[
\binom{tr}{tk}=\sum_{a_1+\dots+a_t=tk}\binom{r}{a_1}\binom{r}{a_2}\dots\binom{r}{a_t}\le(r+1)^{t-1}\binom{r}{k}^{t}.
\]
\end{proof}

\section{Balanced cycles with at least four partite classes}\label{sec:multipartite}
 We begin with some notation and two auxiliary propositions. For a vertex $v$ of a multipartite tournament and a partite class $X$, let $ s(v,X)=d^+(v,X)-d^-(v,X)$. Clearly, $s(v,X)=0$ if $v\in X$.

\begin{proposition}\label{lem:tri-identity}
Let $G$ be a $c$-partite tournament with partite classes
$V_1,\ldots,V_c$, each of size $r$. For any distinct
$i,j,k\in[c]$, the number of triangles in
$G[V_i\cup V_j\cup V_k]$ is 
\[
\frac14\left(
r^3-\sum_{x\in V_i}s(x,V_j)s(x,V_k)
-\sum_{y\in V_j}s(y,V_i)s(y,V_k)
-\sum_{z\in V_k}s(z,V_i)s(z,V_j)
\right).
\]
\end{proposition}

\begin{proof}
Let $t$ denote the number of triangles in
$G[V_i\cup V_j\cup V_k]$.
Every triple $(x,y,z)\in V_i\times V_j\times V_k$
induces either a triangle or a transitive tournament.
Hence
\[
r^3-t
=
\sum_{x\in V_i}d^+(x,V_j)d^+(x,V_k)
+\sum_{y\in V_j}d^+(y,V_i)d^+(y,V_k)
+\sum_{z\in V_k}d^+(z,V_i)d^+(z,V_j).
\]
For every $i\ne j$ with $v\in V_i$, we have $2d^+(v,V_j)=r+s(v,V_j).$ Thus, we obtain
\[
4t
=
r^3-\sum_{x\in V_i}s(x,V_j)s(x,V_k)
-\sum_{y\in V_j}s(y,V_i)s(y,V_k)
-\sum_{z\in V_k}s(z,V_i)s(z,V_j).
\]
\end{proof}

\begin{proposition}\label{lem:dense-cyclic-tuples}
For all integers $c\ge4$ and $k\ge1$
and every $\xi>0$, there exists
$r_0=r_0(c,k,\xi)$ such that the following holds.
Let $G$ be a $c$-partite tournament with
partite classes $V_1,\ldots,V_c$, each of size
$r\ge r_0$. If $G$ contains at least $\xi r^c$ balanced cycles of the form $v_1v_2\cdots v_cv_1$,
where $v_i\in V_i$ for every $i\in[c]$,
then $G$ contains a $k$-balanced cycle.
\end{proposition}

\begin{proof}
Define a $c$-partite $c$-uniform hypergraph $\mathcal{H}$ with vertex classes $V_1,\ldots,V_c$ by declaring that $\{v_1,\ldots,v_c\}\in E(\mathcal H)$ if and only if $v_1\dots v_cv_1$ is a balanced cycle in $G$. We have $e(\mathcal H)\ge\xi r^c.$
By \Cref{thm:erdos-stone-multipartite},
there exist sets
$U_i=\{u_{i,1},\ldots,u_{i,k}\}\subseteq V_i$
for $i\in[c]$ that form a copy of
$K^{(c)}_{k,\ldots,k}$ in $\mathcal H$.
By the definition of $\mathcal H$, every edge between $U_i$ and $U_{i+1}$ is directed from $U_i$ to $U_{i+1}$ for every
$i\in[c]$, where indices are taken modulo $c$.
Consequently,
\[
u_{1,1},u_{2,1},\ldots,u_{c,1},
u_{1,2},u_{2,2},\ldots,u_{c,2},
\ldots,
u_{1,k},u_{2,k},\ldots,u_{c,k},u_{1,1}
\]
forms a $k$-balanced cycle in $G$.
\end{proof}

\subsection{Long cycles}\label{seclarge}
\begin{lemma}\label{large}
There exists an integer $N$ such that
every regular $c$-partite tournament $G$ with
$c\ge4$ and common partite size $r$ contains a
$k$-balanced cycle whenever
$1\le k\le r$ and $ck\ge N$.
\end{lemma}

\begin{proof}
Let $V_1,\ldots,V_c$ be the partite classes of $G$.
For each $i\in[c]$, choose a $k$-subset $U_i$ of $V_i$
uniformly at random.
Set $U=\bigcup_{i=1}^c U_i$ and $H=G[U]$.
For each fixed $v\in V_i$, we have $\mathbb E[d_G^+(v,U)]
=\frac{k}{r}d_G^+(v)=\frac{(c-1)k}{2}.$
Thus, by \cref{cher}, for each fixed $v\in V(G)$
\[
\mathbb P\left(
\left|d_H^+(v)-\frac{(c-1)k}{2}\right|
>\frac{(c-3)k}{4}
\right)
\le
2\exp\left(-\frac{(c-3)^2k}{8(c-1)}\right).
\]
By the union bound, and taking $N$ sufficiently large, we obtain
 $2ck\cdot\exp\left(-\frac{(c-3)^2k}{8(c-1)}\right)\le 2ck\cdot\exp(-\frac{ck}{96})<1$.
Hence there exists $H\subseteq G$ satisfying $\left|d_H^+(v)-\frac{(c-1)k}{2}\right|
\le\frac{(c-3)k}{4}$
for every $v\in V(H)$. Since
$d_H^+(v)+d_H^-(v)=(c-1)k$, we obtain $\dmax(H)-\dsemi(H)
\le\frac{(c-3)k}{2}
<\frac{(c-3)k+2}{2}.$
By \Cref{thm:yeo-hamilton-tool}, $H$ contains a
Hamilton cycle, which is $k$-balanced in $G$.
\end{proof}

\subsection{Short cycles}\label{sec:short-unified}

\begin{lemma}\label{short5}
For all integers $c\ge5$ and $K\ge1$, there exists an integer $r_0=r_0(c,K)\ge K$ such that every
regular $c$-partite tournament with common partite size $r\ge r_0$ contains a $k$-balanced cycle for every
$1\le k\le K$.
\end{lemma}

\begin{proof}
Let $V_1,\ldots,V_c$ be the partite classes of $G$.
Choose one vertex $v_i$ from each class $V_i$ independently and uniformly at random, and let $T$ be the induced tournament.
We claim that $T$ is strongly connected with probability at least $1/100$. 

We first bound the probability that $T$
has a source or a sink. By \cref{reg}, for $v\in V_i$ we have
\[\mathbb P(v\text{ is a source}\mid v\in V(T))
=\prod_{j\ne i}\frac{d^+(v,V_j)}r
\le 2^{-(c-1)}.\] 
The same bound holds for the probability that $v$ is a sink. By the union bound we have
\begin{equation}\label{eq:five-source}
\mathbb P(T\text{ has a source or a sink})
\le\frac{2c}{2^{c-1}}.
\end{equation}

If $c=5$, a tournament with neither a source
nor a sink is strongly connected, since otherwise
its first and last strong components each have
at least three vertices. Hence $\mathbb P(T\text{ is strongly connected})
\ge3/8>1/100.$
Now suppose that $c\ge6$.
If $T$ is not strongly connected and has neither a source nor a sink, let $Q$ be its first strong component. Then $|Q|\ge3$ and $c-|Q|\ge3$. Moreover, $\sum_{v\in Q} d_T^+(v)=e(Q)+|Q|(c-|Q|)$. Consequently,

\begin{equation*}
\resizebox{\linewidth}{!}{$\displaystyle
\begin{aligned}
\sum_{v\in V(T)}\left(d_T^+(v)-\frac{c-1}{2}\right)^2
&= \sum_{v\in Q}\left(d_T^+(v)-\frac{c-1}{2}\right)^2+\sum_{v\notin Q}\left(d_T^+(v)-\frac{c-1}{2}\right)^2\\
&\ge \frac{\left(\sum_{v\in Q}\left(d_T^+(v)-\frac{c-1}{2}\right)\right)^2}{|Q|}+\frac{\left(\sum_{v\notin Q}\left(d_T^+(v)-\frac{c-1}{2}\right)\right)^2}{c-|Q|}=\frac{|Q|(c-|Q|)c}{4}\\
&\ge \frac{3c(c-3)}{4}.
\end{aligned}
$}
\end{equation*}

For each $v\in T$, we have $\mathbb{E}(d_T^+(v))=(c-1)/2$, and $\mathbb{E}((d_T^+(v)-(c-1)/2)^2)=\operatorname{Var}(d_T^+(v))\le (c-1)/4$.
 Hence, by Markov's inequality 
 \begin{equation}\label{nss}
     \mathbb{P}(\sum_{v\in V(T)}\left(d_T^+(v)-\frac{c-1}{2}\right)^2\ge\frac{3c(c-3)}{4})\le \frac{4\mathbb{E}(\sum_{v\in V(T)}\left(d_T^+(v)-\frac{c-1}{2}\right)^2)}{3c(c-3)}\le\frac{c-1}{3(c-3)}.
 \end{equation}

Combining \eqref{eq:five-source} and \eqref{nss} gives $\mathbb P(T\text{ is strongly connected})
\ge1-\frac{2c}{2^{c-1}}-\frac{c-1}{3(c-3)}
\ge\frac1{100}.$
Therefore, at least $r^c/100$ sets induce a strongly connected tournament, and hence contain a  Hamilton cycle by \Cref{thm:camion-tool}.
By the pigeonhole principle, after relabeling the classes if necessary, at least $r^c/(100(c-1)!)$ cycles with the cyclic order $(v_1, v_2,\cdots, v_c,v_1)$ where $v_i\in V_i$. Applying \Cref{lem:dense-cyclic-tuples} with $\xi=1/(100(c-1)!)$ yields the required $k$-balanced cycles.
 \end{proof}

The case $c=4$ differs slightly from the case $c\ge 5$. Volkmann~\cite{vo1999} constructed an infinite family of regular $4$-partite tournaments containing no balanced
directed $4$-cycle.
Thus the preceding lemma cannot be extended to $c=4$ with $k=1$ included.
However, the following lemma shows that $k=1$ is the only exception.

\begin{lemma}\label{short4}
For every integer $K\ge2$, there exists an
integer $r_0=r_0(K)\ge K$ such that every regular $4$-partite tournament with common partite size $r\ge r_0$ contains a $k$-balanced cycle for every
$2\le k\le K$.
\end{lemma}

 \begin{proof}
 Choose $1/r\ll\sigma\ll\eta\ll 1/K\ll1$.
Write $V_1,V_2,V_3,V_4$ for the partite classes, and let $b$ denote the number of balanced $4$-cycles.
If $b\ge\sigma r^4$, the conclusion follows
from \Cref{lem:dense-cyclic-tuples} with parameters $\xi=\sigma/6$ and $c=4$.
Hence we may assume that $b<\sigma r^4$.

\begin{claim}\label{c1}
    We have $\sum_{v\in V(G)}\sum_{i\in[4]}s(v,V_i)^2
    \le 8\sigma r^3.$
    Moreover, for all distinct $i,j\in[4]$, $\sum_{v\in V_i}\bigl|d^\pm(v,V_j)-r/2\bigr|
    \le \sqrt{2\sigma}\,r^2.$
\end{claim}

\begin{proof}
\renewcommand*{\qedsymbol}{$\blacksquare$}
By \Cref{lem:tri-identity}, the number of triangles in $G$ is $r^3+\frac18\sum_{v\in V(G)}
 \sum_{i\in [4]}s(v,V_i)^2$.
On the other hand, every balanced $4$-set induces at most one triangle, unless it contains a $4$-cycle, in which case it induces two. Since each triangle is contained in exactly $r$ balanced $4$-sets, the number of triangles is at most $(r^4+b)/r\le r^3+\sigma r^3$. Thus, $\sum_{v\in V(G)}\sum_{i\in [4]}s(v,V_i)^2\le 8\sigma r^3$.
Moreover,
\begin{equation*}\label{eq:four-degree-error}
\sum_{v\in V_i}\left|d^+(v,V_j)-\frac r2\right|
=\frac12\sum_{v\in V_i}|s(v,V_j)|
\le\frac12
\sqrt{r\sum_{v\in V_i}s(v,V_j)^2}
\le\sqrt{2\sigma}\,r^2.
\end{equation*}
The indegree estimate follows from
$d^+(v,V_j)+d^-(v,V_j)=r$.
\end{proof}

For each $i\in[4]$ and $i\ne j$, let $A_i:=\bigcup_{j}\{v\in V_i:\ |s(v,V_j)|>2\eta r\}$. By \cref{c1}, we have $|A_i|\le \sum_{v\in V(G)}\sum_{i\in [4]}s(v,V_i)^2/(2\eta r)^2\le\eta r$.
For an edge $xy$, let $N_{xy}$ denote the number of balanced $4$-cycles containing $xy$, and let $\mathcal E=\{xy\in E(G):N_{xy}>\eta r^2\}.$
Since $\sum_{xy\in E(G)}N_{xy}=4b$ and $\sigma\ll \eta$, we have
$|\mathcal E|\le4\sigma r^4/\eta r^2\le\eta^2r^2$. 
 Let $B_i$ be the set of vertices in $V_i$ incident with more than $\eta r$ edges of $\mathcal E$. Since $|B_i|\eta r
\le\sum_{v\in V_i}d_{\mathcal E}(v)
\le|\mathcal E|$, we have $|B_i|\le\eta r$.
Let $V_i'=V_i\setminus(A_i\cup B_i)$ and we have for every $v\in V_i'$ and $j\ne i$,
\begin{equation}\label{dpm}
    \left(\frac12-3\eta\right)r\le d^\pm(v,V_j')\le\left(\frac12+\eta\right)r,
\end{equation}
since $|d^\pm(v,V_j)-r/2|=|s(v,V_j)|/2\le\eta r$ for each $v\in V_i$ and at most $2\eta r$ vertices are removed
from $V_j$.

An edge is \emph{good} if both end-vertex lie in $\bigcup_iV_i'$ and the edge does not belong to $\mathcal E$.
A path or cycle is good if all its edges are good.

\begin{claim}\label{c2}
Every good edge $xy$ can be replaced by at least
$r^2/4$ good balanced $(x,y)$-paths of length $3$.
\end{claim}

\begin{proof}
\renewcommand*{\qedsymbol}{$\blacksquare$}
Let $xy$ be a good edge. By relabelling the classes,
we may assume that $x\in V_1'$ and $y\in V_2'$.
Let $\mathcal{P}_{xy}$ be the family of  $(x,y)$-paths
of length $3$ with one internal vertex in each of
$V_3,V_4$, and let $\mathcal{P}_{yx}$ be the family of such
 $(y,x)$-paths.
For $i\in\{3,4\}$, let $X_i=N^+(x)\cap V_i$ and $Y_i=N^+(y)\cap V_i.$ 
By definition,
$|\mathcal{P}_{xy}|=e_G(X_3,V_4\setminus Y_4)+e_G(X_4,V_3\setminus Y_3)$ and $|\mathcal{P}_{yx}|=e_G(Y_3,V_4\setminus X_4)+e_G(Y_4,V_3\setminus X_3)$. Consequently,
\begin{equation*}\label{pq}
|\mathcal{P}_{xy}|+|\mathcal{P}_{yx}|=e_G(X_3,V_4)+e_G(Y_3,V_4)+e_G(X_4,V_3)+e_G(Y_4,V_3)-|X_3||Y_4|-|X_4||Y_3|
\end{equation*}
By \cref{c1}, for every $Z\in\{X_3,Y_3\}$ we have $e_G(Z,V_4)\ge
\frac r2|Z|-\sqrt{2\sigma}\,r^2,$
and similarly, for $Z\in\{X_4,Y_4\}$, $e_G(Z,V_3)\ge
\frac r2|Z|-\sqrt{2\sigma}\,r^2.$ Moreover, since $xy$ is a good edge, we have $\left||X_i|-\frac r2\right|\le\eta r$ and $\left||Y_i|-\frac r2\right|\le\eta r$
for $i\in\{3,4\}$. Therefore, $|\mathcal{P}_{xy}|+|\mathcal{P}_{yx}|$ is at least
\begin{align*}
&\frac r2
\bigl(|X_3|+|X_4|+|Y_3|+|Y_4|\bigr)-|X_3||Y_4|-|X_4||Y_3|
       -4\sqrt{2\sigma}\,r^2\\
&\ge \frac{r^2}{2}
-\left(|X_3|-\frac r2\right)
 \left(|Y_4|-\frac r2\right)
-\left(|X_4|-\frac r2\right)
 \left(|Y_3|-\frac r2\right)
-4\sqrt{2\sigma}\,r^2\\
&\ge
\left(\frac12-2\eta^2-4\sqrt{2\sigma}\right)r^2.
\end{align*}
 Since every path in $\mathcal{P}_{yx}$, together with
$xy$, forms a balanced  $4$-cycle, we have $|\mathcal{P}_{yx}|=N_{xy}\le \eta r^2$. By $\sigma\ll \eta$, we have $|\mathcal{P}_{xy}|\ge\left(1/2-\eta-2\eta^2-4\sqrt{2\sigma}\right)r^2\ge (1/2-2\eta)r^2.$

It remains to exclude paths in $\mathcal{P}_{xy}$ that are not good.
First, since $|V_i\setminus V_i'|\le2\eta r$
for each $i\in[4]$, at most $\bigl(|V_3\setminus V_3'|+|V_4\setminus V_4'|\bigr)r
\le4\eta r^2$ of these paths contain an
internal vertex outside $V_3'\cup V_4'$.
Next, since $x,y\notin B_1\cup B_2$,
each of them is incident with at most
$\eta r$ edges of $\mathcal E$.
Hence, at most $2\eta r^2$ paths have
their first or last edge in $\mathcal E$. 
Finally, since $x$ and $y$ are fixed, each edge
between $V_3$ and $V_4$ can be the second edge
of at most one $(x,y)$-path of length $3$.
Thus, there are at most $|\mathcal E|\le\eta^2r^2$
such paths.
Therefore, the number of good balanced
$(x,y)$-paths is at least $|\mathcal{P}_{xy}|-4\eta r^2-2\eta r^2-\eta^2r^2\ge
\left(1/2-10\eta\right)r^2
\ge r^2/4.$
\end{proof}

Next we find a good $4$-cycle in $G[V_1'\cup V_2']$. 
For $v,v'\in V_1'$, let $p_{vv'}$ denote the number
of $(v,v')$-paths of length $2$ with their internal vertex in $V_2'$. Hence, by \eqref{dpm}, we have
\begin{equation}\label{pvv}
    \sum_{v,v'\in V_1'}p_{vv'}
=\sum_{u\in V_2'}d^+(u,V_1')d^-(u,V_1')
\ge\left(\frac14-4\eta\right)r^3
\end{equation}
and $|p_{vv'}-p_{v'v}|=|d^+(v,V_2')-d^+(v',V_2')|\le4\eta r$.
Consequently, by \eqref{pvv}, the number of $4$-cycles satisfies
\[
\sum_{v,v'\in V_1'}p_{vv'}p_{v'v}/2\ge\sum_{v,v'\in V_1'}p_{vv'}^2/2
-2\eta r\sum_{v,v'\in V_1'}p_{vv'}\ge
\frac{\left(\sum_{v,v'\in V_1'}p_{vv'}\right)^2}{2r^2}
-2\eta r^4\ge
\left(\frac1{32}-3\eta\right)r^4
\]
After excluding at most $|\mathcal{E}|r^2$ cycles with at least one edge in $\mathcal{E}$, we obtain a good $4$-cycle in $G[V_1'\cup V_2']$ with $2$ vertices in $V_1'$ and the remaining vertices in $V_2'$, since $\left(1/32-3\eta\right)r^4-|\mathcal{E}|r^2>0$. 

We now extend a good $4$-cycle in
$G[V_1'\cup V_2']$ to a good $k$-balanced
cycle. Here, to extend an edge $v_iv_j$
between two partite classes means to
replace it with a $(v_i,v_j)$-path of length $3$ whose internal vertices lie in the other two classes, one in each. Starting from a good $4$-cycle in $G[V_1'\cup V_2']$, we extend two
vertex-disjoint edges, regardless of their
orientations, to obtain a good $2$-balanced
cycle containing exactly two edges between
$V_1'$ and $V_2'$ and two between $V_3'$ and
$V_4'$. Now suppose that $2\le j<k$ and we have
a good $j$-balanced cycle with the same
property. Extending an edge between $V_1'$
and $V_2'$ creates a new edge between $V_3'$
and $V_4'$, and vice versa. 
Thus, after extending one edge between $V_1'$ and $V_2'$ and one edge between $V_3'$ and $V_4'$, we obtain a good $(j+1)$-balanced cycle with the same property. By \cref{c2}, these extensions can be
performed greedily while avoiding all
previously used vertices, since at each
step there are at least
$\frac{r^2}{4}-4Kr\ge\frac{r^2}{8}>0$
available replacement paths.
This yields a good $k$-balanced cycle
for every $2\le k\le K$.
\end{proof}

Let $N$ be the constant given by \Cref{large}. There are only finitely many pairs $(c,k)$ satisfying $ck\le N$.
Let $r_0^*$ be the maximum of all the $r_0$ provided
by \Cref{short4,short5} over all such pairs. Combining \Cref{short4,short5,large} completes the proof of \Cref{thm:unified} for $c\ge4$.

\section{Balanced cycles in regular tripartite tournaments}\label{sec:tripartite}

\subsection{The extremal family}\label{vex}

We use the following extremal family of regular tripartite tournaments which are not
robust outexpanders introduced in \cite{DBLPPS}.

\begin{definition}[\cite{DBLPPS}]\label{def:extremal-models}
For $\beta \in [0,1/2]$, $r \in \mathbb{N}$, and disjoint vertex sets $A = \hat{A} \cup \check {A}$, $B$ and $C$ satisfying
$|A| = |B| = |C| = r$,
$|\hat{A}| = (1-\beta)r$ and
$|\check{A}| = \beta r$,
let
$\mathcal{G}_{\beta}(\hat{A},
\check{A}, B, C)$
be the family of tripartite tournaments on
$A\cup B\cup C$ with edge set $E(C,\hat{A})\cup E(\hat{A},B)\cup E(B,\check{A})\cup E(\check{A},C)$
together with a $\beta r$-regular graph oriented
from $C$ to $B$ and its complement
(which has to be $(1-\beta)r$-regular)
oriented from $B$ to $C$.
\end{definition}

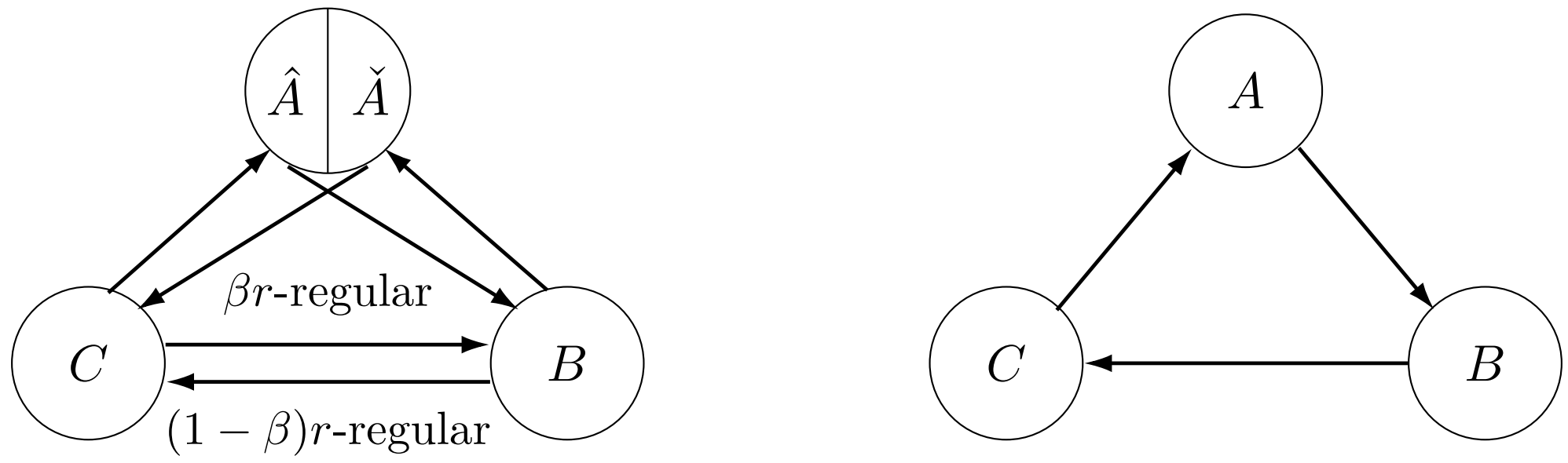
\begin{figure}[htbp]
\centering
\begin{tikzpicture}[
  scale=1.12,
  transform shape,
  >=Latex,
  mainarrow/.style={-{Latex[length=2.4mm,width=1.5mm]},
                    line width=.9pt},
  vertex/.style={circle,draw,minimum size=11.5mm,
                 inner sep=0pt},
  every node/.style={font=\normalsize}
]
\begin{scope}
  \node[vertex] at (-1.8,0) {$C$};
  \node[vertex] at (1.8,0) {$B$};
  \draw (0,2.05) circle[radius=.62];
  \draw (0,1.43) -- (0,2.67);
  \node[font=\normalsize] at (-.32,2.05) {$\hat A$};
\node[font=\normalsize] at (.32,2.05) {$\check A$};

  \draw[mainarrow] (-1.65,.53) -- (-.42,1.62);
  \draw[mainarrow] (-.3,1.48) -- (1.44,.4);
  \draw[mainarrow] (1.65,.55) -- (.42,1.62);
  \draw[mainarrow] (.3,1.48) -- (-1.44,.4);

  \draw[mainarrow] (-1.22,.14) -- (1.22,.14);
  \node[font=\small,fill=white,inner sep=1pt]
    at (0,.52) {$\beta r$-regular};

  \draw[mainarrow] (1.22,-.14) -- (-1.22,-.14);
  \node[font=\small,fill=white,inner sep=1pt]
    at (0,-.52) {$(1-\beta)r$-regular};

  \node at (0,-1.18) {(a) $\mathcal G_\beta$};
\end{scope}

\begin{scope}[xshift=6.9cm]
  \node[vertex] at (0,2.05) {$A$};
  \node[vertex] at (-1.8,0) {$C$};
  \node[vertex] at (1.8,0) {$B$};

  \draw[mainarrow] (-1.42,.40) -- (-.40,1.62);
  \draw[mainarrow] (.40,1.62) -- (1.42,.40);
  \draw[mainarrow] (1.22,0) -- (-1.22,0);

  \node at (0,-1.18) {(b) $G_0$};
\end{scope}
\end{tikzpicture}
\caption{Extremal regular tripartite tournaments: the family $\mathcal G_\beta$
and cyclic extremal tournament $G_0$.}
\label{fig:extremal-models}
\end{figure}

The extremal family $\mathcal G_\beta$ is illustrated
in Figure~\ref{fig:extremal-models}(a).
In particular, when $\beta=0$, it reduces to the
cyclic tripartite tournament $G_0$, with
$C\to A\to B\to C$, as shown in
Figure~\ref{fig:extremal-models}(b).

For two regular tripartite tournaments $G_1$ and $G_2$
with the same partition
$V=V_1\cup V_2\cup V_3$, where $|V_i|=r$ for
each $i\in[3]$, define $\dist(G_1,G_2)
:=\frac12|E(G_1)\triangle E(G_2)|.$
We say that $G_1$ and $G_2$ are $\varepsilon$-close
if $|E(G_1)\triangle E(G_2)|
\le \varepsilon |V|^2.$
Moreover, $G$ is said to be $\varepsilon$-close to
$\mathcal{G}_{\beta}$ if there exists
$G^*\in\mathcal{G}_{\beta}$ that is
$\varepsilon$-close to $G$.
For $\gamma>0$, let $U_\gamma(G_1;G_2)
:=\bigl\{v\in V:
d^+_{G_1\setminus G_2}(v)
+d^-_{G_1\setminus G_2}(v)
\ge \gamma |V|\bigr\}$ and we will often write $U_{\gamma}$ in place of $U_{\gamma}(G_1,G_2)$.
When $G_1$ and $G_2$ are clear from the context, we abbreviate $U_\gamma(G_1;G_2)$ to $U_\gamma$.
The following bound is immediate.

\begin{fact}\label{fact1}
Suppose that $G_1$ and $G_2$ are $\varepsilon$-close. Then $|U_{\gamma}(G_1,G_2)|\le 3\varepsilon r/\gamma$.
\end{fact}
 \begin{proof}
 Since $2|E(G_1)\setminus E(G_2)|=|E(G_1)\triangle E(G_2)|$, we have 
 \[|U_\gamma(G_1;G_2)|\cdot\gamma |V|\le\sum_{v\in U_{\gamma}(G_1,G_2)}d_{G_1\setminus G_2}(v)\le 2|E(G_1)\setminus E(G_2)|\le \varepsilon |V|^2.\] Therefore $|U|\le3\varepsilon r/\gamma$.
 \end{proof}

We use the stability theorem of Di Braccio et
al.~\cite{DBLPPS} in the following form.

\begin{theorem}[Di Braccio et al. \cite{DBLPPS}]
\label{thm:stability}
Let $0<\nu\ll\tau\ll\eps\ll1$.
Then there exists an integer
$r_0$ such that every regular tripartite tournament
with common partite size $r\ge r_0$ is either a robust
$(\nu,\tau)$-outexpander or $\eps$-close to
$\cG_{\beta}$ for some $\beta\in[0,1/2]$.
\end{theorem}

The next two lemmas handle the two alternatives in \Cref{thm:stability}. The robust outexpander case is treated in \Cref{lem:robust-case}, while the case in which $G$ is $\varepsilon$-close to $\mathcal G_\beta$ is treated in \Cref{lem:extremal-transfer}. Together, they establish \Cref{thm:unified} for $c=3$.

\subsection{Linear forests}\label{lfor}

For later use, we combine weaker versions of Lemmas~4.11
and~4.12 of Di Braccio et al.~\cite{DBLPPS} into a single lemma. Those results provide linear many edge-disjoint forests, here we need only one.

\begin{lemma}[Di Braccio et al.~\cite{DBLPPS}]
\label{forest}
Let $0<1/r\ll\varepsilon\ll\gamma\ll1$.
Suppose that $G$ is a regular tripartite tournament
with partite classes $A,B,C$, each of size $r$,
and that $G$ is $\varepsilon$-close to
$G'\in\mathcal{G}_{\beta}(\hat A,\check A,B,C)$,
where $A=\hat A\cup\check A$ and either
$\beta=0$ or $\gamma\ll\beta\le1/2$.
Then $G$ contains a linear forest $F$ such that $e(F)\le\gamma r$ and $d_F^+(v)=d_F^-(v)=1$ for every $v\in U_\gamma(G,G')$.
Moreover, the following hold:

\textnormal{($i$)} If $\beta=0$, then $e_F(C,B)=e_F(B,A)=e_F(A,C).$

\textnormal{($ii$)} If $\beta>0$, then $e_F(C,A)=e_F(A,B)$ and $e_F(B,A)=e_F(A,C)$.
Moreover, each path in $F$ is a $(C,B)$-path.
\end{lemma}

We also use two results on Hamilton cycles from the same paper.

\begin{lemma}[Di Braccio et al. \cite{DBLPPS}]\label{hcc2}
Let $1/n\ll\eta\ll1$.  Let $G$ be a tripartite digraph with partite classes $A,B,C$ such that $(1-\eta)n\le|A|\le n$ and $|B|=|C|=n$.  Suppose every
vertex has at least $(1-\eta)n$ outneighbors in the next part and at least
$(1-\eta)n$ inneighbors in the previous part with respect to
$C\to A\to B\to C$.  Then every matching $M\subseteq E_{G}(C,B)$ of size exactly
$n-|A|$ is contained in a Hamilton cycle of $G$.
\end{lemma}

\begin{lemma}[Di Braccio et al. \cite{DBLPPS}]\label{hcc1}
Let $1/n\ll\eta\ll1$.  Let $G$ be a four-partite digraph with partition
$V(G)=X\cup Y\cup B\cup C$, where $|B|=|C|=n$ and
$|X|,|Y|\le(1-8\eta)n$.  Suppose
\begin{align*}
 \dplus_G(x,B),\dminus_G(x,C)&\ge(1-\eta)n &&(x\in X),\\
 \dplus_G(y,C),\dminus_G(y,B)&\ge(1-\eta)n &&(y\in Y),\\
 \dplus_G(b,Y\cup C),\dminus_G(b,X\cup C)&\ge(1-\eta)n &&(b\in B),\\
 \dplus_G(c,X\cup B),\dminus_G(c,Y\cup B)&\ge(1-\eta)n &&(c\in C).
\end{align*}
Then every matching $M\subseteq E_{G}(C,B)$ with
$|M|\le\eta n$ is contained in a Hamilton cycle of $G$.
\end{lemma}

\subsection{Non-extremal case}\label{nex}

\begin{lemma}\label{lem:robust-case}
Let $0<\nu\le\tau\ll1$.
There exists an integer $r_0$ such that every regular tripartite
tournament $G$ with partite classes $A,B,C$ of common size
$r\ge r_0$ which is a robust $(\nu,\tau)$-outexpander contains
a $k$-balanced cycle for every integer $1\le k\le r$.
\end{lemma}

\begin{proof}
Let $A=V_1,B=V_2,C=V_3$. Choose constants satisfying
\[
 0<1/r_0\ll\xi\ll\nu'\ll1/M\ll1/M'\ll\varepsilon
 \ll d,\nu,\tau\ll 1.
\]

\noindent \textbf{Case 1.} $k\le\xi r$

By \cref{lem:tri-identity}, $G$ contains at least
$r^3/4$ triangles. Relabeling the classes if necessary, we may assume that at least $r^3/8$ of them have orientation $A\to B\to C\to A$.
Let $H$ be the undirected tripartite graph with $V(H)=V(G)$ and $E(H)=E_{G}(A,B)\cup E_{G}(B,C)\cup E_{G}(C,A)$. Thus $H$ contains at least $r^3/8$ triangles.

We apply \cref{lem:regularity} to $H$ with parameters $(\varepsilon,d,M')$ and partite classes $A,B,C$, obtaining a spanning subgraph $H'$, an exceptional set $V_0$ with $|V_0|\le\varepsilon n$, and $3t$ clusters $A_1,\ldots,A_t,B_1,\dots,B_t,C_1,\dots,C_t$ of equal size $m$, where $M'\le t\le M$.
Let $R$ be the corresponding reduced graph with parameters $(\varepsilon,d)$.
 By \cref{lem:regularity} (iii), $e(H)-e(H')=\frac12\sum_{x\in V(H)}
       \bigl(d_{H}(x)-d_{H'}(x)\bigr)
 \le\frac12(d+\varepsilon)n^2.$
Each deleted edge belongs to at most $r$ triangles, and each exceptional vertex belongs to at most $r^2$ triangles.
Hence $H'-V_0$ contains at least
\[
 \frac{r^3}{8}
 -\frac12(d+\varepsilon)n^2r
 -\varepsilon nr^2
 =
 \left(\frac18-\frac92d-\frac{15}{2}\varepsilon\right)r^3
 >0
\]
triangles.
Therefore $R$ contains a triangle $T$. We may assume that its vertices are $V_{11},V_{21},V_{31}$.

Applying \cref{prop:super-regularization} to the triangle $T$ with $\Delta = 2$, we obtain a subgraph $H_{T}$ of $H'$.
More precisely, $V(H_T)=V_{11}'\cup V_{21}'\cup V_{31}'$, where each $V_{i1}'\subseteq V_{i1}$ is obtained by deleting exactly
$2\varepsilon m$ vertices. As $t\le M$ and $tm\ge(1-\varepsilon)r$, we have $|V_{i1}'|\ge(1-2\varepsilon)(1-\varepsilon)r/M\ge\xi r\ge k.$ Furthermore, the three unordered pairs $(V_{i1}',V_{j1}')$ are all
$(\sqrt{\varepsilon},d-3\varepsilon)$-super-regular.
By \cref{blowup}, one can embed a $k$-balanced cycle into $H'$. By the definition of $H$, this corresponds to a $k$-balanced cycle in $G$.

\noindent \textbf{Case 2.} $k>\xi r$

Independently for each $i\in[3]$, choose a $k$-subset
$U_i$ of $V_i$ uniformly at random, and let
$U=U_1\cup U_2\cup U_3$.
Consider the $3k$-subsets $S$ of $V(G)$ for which $G[S]$ is not a robust $(\nu',4\tau)$-outexpander.
Since $\delta^0(G)=r=|V(G)|/3$,
\cref{lem:random-robust-expansion,log} imply that
the number of such sets $S$ is at most
\[
\frac{1}{9r^2}\binom{3r}{3k}\le\frac{1}{9r^2}\sum_{\substack{0\le a,b,c\le r\\a+b+c=3k}}
  \binom{r}{a}\binom{r}{b}\binom{r}{c}\le\frac{(r+1)^2}{9r^2}\binom{r}{k}^{3}
\]
Since there are exactly $\binom{r}{k}^3$
choices for $U$, it follows that
\[
\mathbb P\bigl(
G[U]\text{ is not a robust }
(\nu',4\tau)\text{-outexpander}
\bigr)
\le\frac{(r+1)^2}{9r^2}
\le\frac14.
\] 
For each fixed $v\in V_i$, we have $\mathbb E[d_G^\pm(v,U)]
=\frac{k}{r}\sum_{j\in[3]\setminus\{i\}}
d_G^\pm(v,V_j)=k.$
By \cref{cher} and the union bound, $\mathbb P\bigl(\dsemi(G[U])<3k/4\bigr)
 \le4\cdot 3r e^{-k/16}< 1/4$.
Consequently, there exists a $3k$-set $U$ such that
$G[U]$ is a robust $(\nu',4\tau)$-outexpander and $\dsemi(G[U])\ge\frac{3k}{4}
=\frac{|U|}{4}.$

Since $\nu'\le4\tau<1/64$,
\cref{thm:lo-patel} gives a Hamilton cycle in $G[U]$. This cycle contains exactly $k$ vertices from each $V_i$ and is therefore $k$-balanced.
\end{proof}

\subsection{The extremal case}\label{ex}

\begin{lemma}\label{lem:extremal-transfer}
There exists $\eps>0$ such that, for all sufficiently large $r$, every regular
tripartite tournament $G$ which is $\eps$-close to $\cG_{\beta}$ for some
$0\le \beta\le 1/2$ contains a $k$-balanced cycle for every $1\le k\le r$.
\end{lemma}

\begin{proof}
Choose constants satisfying
\begin{equation}\label{p1}
    1/r\ll
 \eps\ll\gamma\ll\lambda\ll\eta
 \ll\sigma\ll\tau\ll1.
\end{equation}
Let $G_0$ be the cyclic extremal tournament with $A\to B\to C\to A$. Choose $G_1\in\mathcal G_\beta(\hat A,\check A,B,C)$ such that $\operatorname{dist}(G,G_1)\le 9\varepsilon r^2/2$. Let $d=\beta r$ and $q=(1-\beta)r$.

\medskip
\noindent \textbf{Case 1. }$k\le (1-\lambda) r.$

Let $U=U_\gamma(G,G_1)$. For each $V\in\{\hat A,\check A,B,C\}$, write $V_0=V\setminus U$.
By \eqref{p1} and \cref{fact1} for every $b\in B_0$ and $c\in C_0$, the following hold 
\begin{equation}\label{b1}
d_G^+(b,C_0),d_G^-(c,B_0)\ge q-3\gamma r-|U|
\ge q-\lambda r/20>0\qquad\text{and}
\end{equation}
\begin{equation}\label{b2}
d_G^-(b,C_0),d_G^+(c,B_0)\ge d-3\gamma r-|U|
\ge d-\lambda r/20.
\end{equation}

If $q-\lambda r/20<k$, then $d-\lambda r/20>0$.
By \eqref{b1} and \eqref{b2}, there exist $c\in C_0$ and $b_1,b_2\in B_0$ such that $b_1c,cb_2\in E(G)$. Consider disjoint subsets $B_1,B_2$ of $B_0$ containing $b_1,b_2$, respectively, with $|B_1|=q-\lambda r/20$ and $|B_2|=k-q+\lambda r/20\le d-\lambda r/20$. By Hall's theorem, we obtain matchings $M_1=M[B_1,C_1]$ and $M_2=M[C_2,B_2]$ containing $b_1c,cb_2$, respectively.
By construction, the matchings $M_1$ and $M_2$ both contain the vertex $c$. Let
$C_1\cap C_2=\{v_1,\ldots,v_t\}$. Since $|C_0|\ge k$ and $|C_1\cup C_2|=k-t$, we may choose distinct
vertices $u_1,\ldots,u_t\in C_0\setminus(C_1\cup C_2)$.
Distribute the vertices of $C_2\setminus C_1$ arbitrarily between the pairs $v_i,u_i$, and those of $C_1\setminus C_2$ between the pairs $u_i,v_{i+1}$, where $v_{t+1}=v_1$. The number and order of vertices placed between each pair are arbitrary. This yields a cyclic ordering $\pi$ of the $k$ vertices of $C_0$.
Call an ordered pair $(c,c')$ of consecutive vertices in this cyclic ordering type I if $c\in C_2$, and of type II if $c'\in C_1$.
By construction, no consecutive ordered pair $(c,c')$
in this cyclic ordering $\pi$ satisfies $c\in C_2$ and $c'\in C_1$. Thus, the two types are disjoint.
Since there are exactly $|C_2|=|B_2|$ type I pairs
and $|C_1|=|B_1|$ type II pairs, and
$|B_1|+|B_2|=k$, they account for all $k$ consecutive ordered pairs.

By the definition of $G_1$, for every
$b\in B_0$ and $c\in C_0$, we have
\begin{equation}\label{hall1}
\begin{aligned}
|\{y\in \check{A}_0:b\to y\to c\}|
&\ge d-|U|-6\gamma r\ge d-\frac{\lambda r}{20}\ge |B_2|,\qquad\text{if} \quad B_2\neq\emptyset
\\
|\{x\in \hat{A}_0:c\to x\to b\}|
&\ge q-|U|-6\gamma r\ge q-\frac{\lambda r}{20}\ge |B_1|.
\end{aligned}
\end{equation}
For each type I pair $(c_1,c_2)$, let $b$ be the vertex matched to $c_1$ by $M_2$.
By \eqref{hall1} and Hall's theorem, we can
choose distinct vertices $y\in \check{A}_0$ to obtain
paths $c_1byc_2$.
Similarly, for each type II pair $(c_3,c_4)$, let $b$ be the vertex matched to $c_4$ by $M_1$. We then choose distinct vertices $x\in\hat A_0$ so as to obtain the paths $c_3xbc_4$.
Replacing all consecutive pairs by these paths
gives a $k$-balanced cycle.

If $q-\lambda r/20\ge k$. Choose a $k$-subset $B_1\subseteq B_0$. As above, there is a matching $M_1=M[B_1,C_1]$ of size $k$. Let $B_2=C_2=\emptyset$ and choose an arbitrary cyclic ordering of $C_1$.
All consecutive pairs are of type II, the same argument as above allows us to replace these pairs with paths, thereby obtaining a $k$-balanced cycle.

\medskip
\noindent \textbf{Case 2.} $(1-\lambda) r\le k\le r$.

We divide this case into two subcases according to whether $\beta<\sigma$ or $\beta\ge\sigma$. Let $m=\lambda r$.

\medskip
\noindent \textbf{Subcase 2.1.} $\beta<\sigma$.

Since $G_1$ differs from $G_0$ in exactly
$3\beta r^2$ edges, we have
$\dist(G,G_0)\le \dist(G,G_1)+\dist(G_1,G_0)\le(\frac92\eps+3\beta)r^2\le\frac92\sigma r^2$.
Thus $G$ is $\sigma$-close to $G_0$. Since $\sigma\ll\tau$,
\Cref{forest}(i) gives a linear forest $F_0$
with $e(F_0)\le\tau r$,
$e_{F_0}(C,B)=e_{F_0}(B,A)=e_{F_0}(A,C)$ and
\begin{equation}\label{UF}
    U_{\tau}(G,G_0)\subseteq \mathring V(F_0)
\end{equation}
Indeed, we can extend each component of $F_0$ to a $(C,B)$-path by adding at most two edges at each end, using edges from $E_G(A,B)$, $E_G(B,C)$, or $E_G(C,A)$. By \eqref{UF} there are at least $(1-3\tau)r$ choices of each added edge, greater than $|V(F_0)|+4e(F_0)\le6\tau r$, so the extensions can be made greedily and vertex-disjointly. Let $F$ be the resulting forest. Then $|V(F)|\le6\tau r$. Moreover, the way in which $F$ was obtained from $F_0$ gives
\begin{equation}\label{ef}
    e_{F}(C,B)=e_{F}(B,A)=e_{F}(A,C).
\end{equation}

Next, we construct a path $R=c_0x_1b_1c_1\cdots c_mx_{m+1}b_{m+1}$ that is vertex-disjoint from $F$. Its vertices satisfy $c_i\in C$,
$x_i\in A$ and $b_i\in B$, and it has the additional property that $c_0x_i\in E(G)$ for every $i\in[m+1]$.
Indeed, by \eqref{UF}, for every
$x\in A\setminus V(F)$, $b\in B\setminus V(F)$ and $c,c'\in C\setminus V(F)$ we have
\[
|N_G^+(c,A)\cap N_G^+(c',A)|\ge r-6\tau r,
\qquad d_G^+(x,B),\ d_G^+(b,C)\ge r-3\tau r.
\]
Starting with $c_0\in C\setminus V(F)$, we successively
choose $x_i\in N_G^+(c_0,A)\cap N_G^+(c_{i-1},A)$,
$b_i\in N_G^+(x_i,B)$, and, if $i\le m$,
$c_i\in N_G^+(b_i,C)$, avoiding used vertices.
This is possible since at most $|V(F)|+3(m+1)\le 7\tau r<(1-6\tau)r$ vertices must be avoided.

Let $\mathcal P=F\cup R$, and write $|\mathcal P|$ for its number of components. By construction, every component of $\mathcal P$ is a $(C,B)$-path.
By the construction of $\mathcal P$ and \eqref{ef}, we have $|\mathcal{P}|\le e(F_0)+1\le \tau r+1$ and $|V(\mathcal{P})|\le7\tau r$.
\begin{equation}\label{P1}
    e_{\mathcal P}(C,B)=e_{\mathcal P}(B,A)=e_{\mathcal P}(A,C).
\end{equation}
Let $G^*=G\circ\mathcal{P}$ be the auxiliary tripartite digraph, whose partite classes $A',B',C'$ are subsets of $A,B,C$, respectively. Let $M_{\mathcal P}$ denote the auxiliary matching, which has size $|\mathcal P|$.
Since every vertex of $V(\mathcal P)\cap A$ lies in
$\mathring V(\mathcal P)$, we have
\[
|V(\mathcal P)\cap A|
=e_{\mathcal P}(C,A)+e_{\mathcal P}(B,A)
=e_{\mathcal P}(A,B)+e_{\mathcal P}(A,C).
\]
Similarly,
$|V(\mathcal P)\cap B|
=e_{\mathcal P}(A,B)+e_{\mathcal P}(C,B)$
and
$|V(\mathcal P)\cap C|
=e_{\mathcal P}(C,A)+e_{\mathcal P}(C,B)$.
Combining these with \eqref{P1}, we obtain
$|V(\mathcal P)\cap A|
=|V(\mathcal P)\cap B|
=|V(\mathcal P)\cap C|=:p$.
Consequently,
$|A'|=r-p$ and $|B'|=|C'|=:r'=r-p+|\mathcal P|$. Since $p\le |V(\mathcal P)|\le7\tau r$, we have
\begin{equation}\label{r''}
    r'\ge r- p\ge(1-7\tau)r.
\end{equation}

We now verify the hypotheses of \Cref{hcc2}
for $G^*$. Since every vertex of $G^*$
lies outside $U_{\tau}(G,G_0)$, for $v\in A'$
\[
d_{G^*}^+(v,B')\ge d_{G}^+(v,B)-(p-|\mathcal{P}|)\ge r-3\tau r-p+|\mathcal{P}|\ge(1-10\tau)r',
\]
where the last inequality follows from
\eqref{p1} and \eqref{r''}.
All other required degree conditions
follow analogously.
Moreover, $(1-10\tau)r'
\le|A'|\le r'$ and
$|M_{\mathcal P}|=|\mathcal{P}|=r'-|A'|$.
Applying \Cref{hcc2} to $G^*$ with parameter $10\tau$, we obtain a Hamilton cycle containing $M_{\mathcal P}$. Replacing each edge of $M_{\mathcal P}$ with the path in $\mathcal P$ yields a Hamilton cycle of $G$ containing $R$.
Finally, for every integer $0\le t\le m$,
replace the subpath $c_0x_1b_1c_1\cdots c_tx_{t+1}$
of $R$ by the edge $c_0x_{t+1}$. This removes exactly $t$ vertices from each of $A,B,C$, yielding an $(r-t)$-balanced cycle. Therefore,
we obtain a $k$-balanced cycle for every $r-\lambda r\le k\le r$.

\medskip
\noindent\textbf{Subcase 2.2.} $\beta\ge\sigma$.

By \Cref{forest}(ii), there exists a linear forest
$F$ whose components are $(C,B)$-paths,
such that $e(F)\le\gamma r$,
$U_\gamma(G,G_1)\subseteq\mathring V(F)$, and
$e_F(C,A)=e_F(A,B)$,
$e_F(B,A)=e_F(A,C)$.
By the same argument as in Subcase 2.1,
we construct a path
$Q=c_0x_1b_1c_1\cdots c_mx_{m+1}b_{m+1}$ 
disjoint from $F$, where $c_i\in C$,
$x_i\in \hat{A}$ and $b_i\in B$, such that
$c_0x_i\in E(G)$ for every $i\in[m+1]$.
Indeed, for any $c,c'\in C\setminus V(F)$,
$x\in\hat A\setminus V(F)$ and
$b\in B\setminus V(F)$, we have $|N_G^+(c,\hat A)\cap N_G^+(c',\hat A)|\ge(1-\beta-6\gamma)r$, $d_G^+(x,B)\ge(1-3\gamma)r$ and $d_G^+(b,C)\ge(1-\beta-3\gamma)r$.
Since $\beta\le 1/2$, by \eqref{p1}, $(1-\beta-6\gamma)r\ge(2\gamma+4\lambda)r\ge |V(F)|+3(m+1).$
Thus, the same greedy construction as in Subcase 2.1
applies.

Let $\mathcal L=F\cup Q$. Then $\mathcal L$ is a linear forest
whose components are $(C,B)$-paths.
Moreover, $U_\gamma(G,G_1)\subseteq\mathring V(\mathcal L)$
and $e_{\mathcal L}(C,A)=e_{\mathcal L}(A,B)$.
We also have $|\mathcal L|\le e(F)+1\le\gamma r+1$ and $|V(\mathcal L)|\le(2\gamma+4\lambda)r$.
Let $G^*=G\circ\mathcal L$, and let $M_{\mathcal L}$
be the auxiliary matching of size $|\mathcal L|$.
As before, we verify the hypotheses of \Cref{hcc1}
for $G^*$. We have
\[
|V(\mathcal L)\cap B|
=e_{\mathcal L}(A,B)+e_{\mathcal L}(C,B)=e_{\mathcal L}(C,A)+e_{\mathcal L}(C,B)
=|V(\mathcal L)\cap C|.
\]
Let $\ell$ denote this common value.
Then $|B'|=|C'|=:r'=r-\ell+|\mathcal L|$. 
Since $\ell\le |V(\mathcal L)|\le(2\gamma+4\lambda)r$, we have
\begin{equation}\label{r'}
    r'\ge r-\ell\ge(1-2\gamma-4\lambda)r.
\end{equation}
Moreover, by \eqref{p1} and \eqref{r'},
$|\hat A'|\le(1-\beta)r\le(1-\sigma)r
\le(1-8\eta)r'$ and
$|\check A'|\le\beta r\le r/2\le(1-8\eta)r'$.
Also, $|M_{\mathcal L}|=|\mathcal L|\le\eta r'$.
Since every vertex of $G^*$ lies outside
$U_\gamma(G,G_1)$, the definition of $G_1$ gives
\[
\begin{aligned}
d_{G^*}^+(v,B'),\ d_{G^*}^-(v,C')
&\ge r'-3\gamma r &&(v\in\hat A'),\\
d_{G^*}^-(v,B'),\ d_{G^*}^+(v,C')
&\ge r'-3\gamma r &&(v\in\check A'),\\
d_{G^*}^+(v,\check A'\cup C'),\
d_{G^*}^-(v,\hat A'\cup C')
&\ge r'-3\gamma r-|V(\mathcal L)\cap A|
&&(v\in B'),\\
d_{G^*}^+(v,\hat A'\cup B'),\
d_{G^*}^-(v,\check A'\cup B')
&\ge r'-3\gamma r-|V(\mathcal L)\cap A|
&&(v\in C').
\end{aligned}
\]
All these bounds are at least $(1-\eta)r'$,
since \eqref{p1} and \eqref{r'}
imply that $3\gamma r+|V(\mathcal L)\cap A|
\le(5\gamma+4\lambda)r
\le\eta r'.$
Thus all the hypotheses of \Cref{hcc1}
are satisfied, and $G^*$ has a Hamilton
cycle containing $M_{\mathcal L}$.
After replacing edges of $M_{\mathcal{L}}$ by $\mathcal{L}$ gives a
Hamilton cycle of $G$ containing $Q$.
As in Subcase 2.1, shortcutting $Q$
yields a $k$-balanced cycle for every
integer $(1-\lambda)r\le k\le r$.
\end{proof}

\section{\texorpdfstring{$(p,q)$}{(p,q)}-Cycle factors}\label{sec:complementary}

\begin{proof}[Proof of \Cref{thm:complementary}]
Let $n=cr$ and $s=\min\{p,n-p\}$. 
Suppose first that $s<n/12$. By
\Cref{thm:yeo-vertex-tool,thm:yeo-four-tool}, $G$ has a cycle $C$ of length $s$. Let $H=G-V(C)$. For every $v\in V(H)$, $|d^{\pm}_{H}(v)-(c-1)r/2|\le s$, so $\dmax(H)-\dsemi(H)\le s$. Let $k_{1}\le k_2$ be the sizes
of the two largest partite classes of $H$. Since
$k_{1},k_2\le r$,
\[
 \frac{|V(H)|-k_1-2k_2+2}{2}
 \ge\frac{(c-3)r-s+2}{2}\ge s,
\]
where the last inequality follows from
$3s<n/4\le(c-3)r$. By \Cref{thm:yeo-hamilton-tool},
$H$ has a Hamilton cycle, together with $C$ forms a $(p,n-p)$-cycle factor.

Now let $s\ge n/12$. Choose a uniformly random $p$-set $U_1\subseteq V(G)$, and let $U_2=V(G)\setminus U_1$. For every fixed $v\in V(G)$, we have $\mathbb E[d_G^\pm(v,U_i)]
=\frac{|U_i|}{n}d_G^\pm(v)
=\frac{|U_i|(c-1)}{2c}
=:\mu_i.$
By \cref{cher}, for every fixed $v\in V(G)$
and $i\in[2]$,
\begin{equation}\label{mu1}
 \mathbb P\bigl(
     |d_G^\pm(v,U_i)-\mu_i|>\frac{|U_i|}{28}\bigr)
 \le 2\exp(-\frac{|U_i|}{400}).
\end{equation}
For each partite class $V_j$, $\mathbb E(|V_j\cap U_i|)=|U_i|/c$, so
\begin{equation}\label{mu2}
 \mathbb P\bigl(\bigl||V_j\cap U_i|-\frac{|U_i|}{c}\bigr|>\frac{|U_i|}{28}\bigr)
 \le 2\exp(-\frac{|U_i|}{400}).
\end{equation}
By the union bound and \eqref{mu1} and \eqref{mu2} there exists a partition $\{U_1,U_2\}$ such that, for each $v\in U_i$, $|d^\pm(v,U_i)-\mu_i|\le |U_i|/28$,
 and $|V_j\cap U_i|\le |U_i|/c+|U_i|/28
 \le2|U_i|/7$ for $i\in [2]$ and $j\in [c]$, because $(8n+4c)\exp(-s/400)<1$ and $c\ge 4$.
Let $m_{i,1}\le m_{i,2}$ be the sizes
of the two largest partite classes of $G[U_i]$.
\[\dmax(G[U_i])-\dsemi(G[U_i])\le \frac{|U_i|}{14}
 \le \frac{|U_i|-m_{i,1}-2m_{i,2}+2}{2}.\]
By \Cref{thm:yeo-hamilton-tool}, each of $G[U_1]$ and
$G[U_2]$ has a Hamilton cycle. These two vertex-disjoint Hamilton cycles form the required $(p,n-p)$-cycle-factor, completing the proof.
\end{proof}

\vspace{0.4cm}
\noindent \textbf{Declaration of AI.}
The authors used GPT-5.6 Sol as a research assistance
tool during the preparation of this manuscript.
In particular, while investigating the case of
regular tripartite tournaments, the AI helped identify an
important reference. The arguments presented in \cref{vex} and \cref{lfor} make substantial use of
techniques and lemmas developed in that reference,
which is acknowledged and cited in the manuscript.
The authors also used AI to assist with language
editing, typesetting, and checking the correctness
of mathematical arguments. The authors have
carefully reviewed the manuscript and take full
responsibility for its mathematical correctness,
accuracy, and content.
\medskip

\noindent \textbf{Acknowledgments.} We thank Yangyang Cheng for bringing these problems to our attention and sharing his insights with us.


{\small
\begin{thebibliography}{99}

\bibitem{BessyThiebaut2023}
S. Bessy and J. Thiebaut,
Complementary cycles of any length in regular bipartite tournaments,
\emph{J. Graph Theory} 103 (2023), 186--211.

\bibitem{Bondy}
J. A. Bondy,
Diconnected orientations and a conjecture of Las Vergnas,
\emph{J. Lond. Math. Soc.} 14 (1976), 277--282.

\bibitem{Camion}
P. Camion,
Chemins et circuits hamiltoniens des graphes complets,
\emph{C. R. Acad. Sci. Paris} 249 (1959), 2151--2152.

\bibitem{DBLPPS}
F. Di Braccio, J. Lada, V. Patel, Y. Pehova, and J. Skokan,
Hamilton decompositions of regular tripartite tournaments,
arXiv:2507.12460 (2025).

\bibitem{E}
P. Erd\H{o}s,
On extremal problems of graphs and generalized graphs,
\emph{Israel J. Math.} 2 (1964), 183--190.

\bibitem{GGKLO}
A. Gir\~ao, B. Granet, D. K\"uhn, A. Lo, and D. Osthus,
Path decompositions of tournaments,
\emph{Proc. Lond. Math. Soc.} 126 (2023), 429--517.

\bibitem{Hoe}
W. Hoeffding, Probability inequalities for sums of bounded random variables,
\emph{J. Amer. Statist. Assoc.} 58 (1963), 13--30.

\bibitem{KSS}
J. Koml\'os, G. N. S\'ark\"ozy, and E. Szemer\'edi,
Blow-up lemma,
\emph{Combinatorica} 17 (1997), 109--123.

\bibitem{KOT2005}
D. K\"uhn, D. Osthus, and A. Taraz,
Large planar subgraphs in dense graphs,
\emph{J. Combin. Theory Ser. B} 95 (2005), 263--282.

\bibitem{LoPatel}
A. Lo and V. Patel,
Hamilton cycles in sparse robustly expanding digraphs,
\emph{Electron. J. Combin.} 25 (2018), Paper 3.44.

\bibitem{Moon}
J. W. Moon,
On subtournaments of a tournament,
\emph{Canad. Math. Bull.} 9 (1966), 297--301.

\bibitem{Redei}
L. R\'edei,
Ein kombinatorischer Satz,
\emph{Acta Litt. Sci. Szeged} 7 (1934), 39--43.

\bibitem{StellaVolkmannNine}
I. Stella and L. Volkmann,
Each regular 3-partite tournament has a 9-cycle,
\emph{An. Univ. Bucure\c{s}ti Mat.} 56 (2007), 91--120.

\bibitem{vo1999}
L. Volkmann,
Strong subtournaments of multipartite tournaments,
\emph{Australas. J. Combin.} 20 (1999), 189--196.

\bibitem{VolkmannFour2004}
L. Volkmann,
Complementary cycles in regular multipartite tournaments, where one cycle has length four,
\emph{Kyungpook Math. J.} 44 (2004), 219--247.

\bibitem{VolkmannComplementary}
L. Volkmann,
All regular multipartite tournaments that are cycle complementary,
\emph{Discrete Math.} 281 (2004), 255--266.

\bibitem{Volkmann}
L. Volkmann,
On cycles in regular 3-partite tournaments,
\emph{Discrete Math.} 306 (2006), 1198--1206.

\bibitem{vp}
L. Volkmann and S. Winzen,
Paths with a given number of vertices from each partite set in regular multipartite tournaments,
\emph{Discrete Math.} 306 (2006), 2724--2732.

\bibitem{VolkmannWinzenCycles}
L. Volkmann and S. Winzen,
Cycles with a given number of vertices from each partite set in regular multipartite tournaments,
\emph{Czechoslovak Math. J.} 56 (2006), 827--843.

\bibitem{VolkmannWinzenStrong}
L. Volkmann and S. Winzen,
Almost regular $c$-partite tournaments contain a strong subtournament of order $c$ when $c\ge5$,
\emph{Discrete Math.} 308 (2008), 1710--1721.

\bibitem{YeoHamilton}
A. Yeo,
One-diregular subgraphs in semicomplete multipartite digraphs,
\emph{J. Graph Theory} 24 (1997), 175--185.

\bibitem{YeoNear}
A. Yeo,
How close to regular must a semicomplete multipartite digraph be to secure hamiltonicity?,
\emph{Graphs Combin.} 15 (1999), 481--493.

\bibitem{YeoVertex}
A. Yeo,
Diregular $c$-partite tournaments are vertex-pancyclic when $c\ge5$,
\emph{J. Graph Theory} 32 (1999), 137--152.

\bibitem{Yeo2007}
A. Yeo,
Paths and cycles containing given arcs, in close to regular multipartite tournaments,
\emph{J. Combin. Theory Ser. B} 97 (2007), 949--963.

\end{thebibliography}
}
\end{document}